\documentclass[12pt]{amsart}

\usepackage[T1]{fontenc}
\usepackage{lmodern}
\usepackage{microtype}
\usepackage{mathtools}
\usepackage{amssymb}
\usepackage{xcolor}
\usepackage{enumitem,tikz}
\usepackage[hidelinks]{hyperref}
\usepackage{comment} 
\usepackage[margin=1.2in]{geometry}

\newtheorem{theorem}{Theorem}[section]
\newtheorem{proposition}[theorem]{Proposition}
\newtheorem{lemma}[theorem]{Lemma}
\newtheorem{corollary}[theorem]{Corollary}

\newtheorem{example}[theorem]{Example}
\newtheorem{definition}{Definition}[section]
\newtheorem{conjecture}[theorem]{Conjecture}

\newcounter{lthm}
\theoremstyle{plain}
\newtheorem{letterthm}[lthm]{Theorem}

\newcommand{\wdim}{\operatorname{w.dim}}
\newcommand{\fd}{\operatorname{fd}}
\newcommand{\Tor}{\operatorname{Tor}}
\newcommand{\Frac}{\operatorname{Frac}}
\newcommand{\gdim}{\operatorname{gl.dim}}

\title[Weak Hilbert--syzygy equality]{Coherence and the Hilbert--syzygy identity for power series over valuation domains}

\author{Mohamed Amine Aouichaoui}
\address{Department of Mathematics, Faculty of Sciences of Monastir, University of Monastir, Monastir, Tunisia}
\email{aouichaoui.mohamedamine@fsm.u-monastir.tn}

\author{Adam Jones}
\address{Department of Pure Mathematics and Mathematical Statistics, Centre for Mathematical Sciences, University of Cambridge, Cambridge United Kingdom}
\email{aj718@cam.ac.uk}

\author{Yiding Wang}
\address{Qiuzhen College, Tsinghua University, Beijing, China}
\email{wyd23@mails.tsinghua.edu.cn}


\begin{document}

\begin{abstract}
Let $V$ be a complete valuation domain of rank 1, and let $\Gamma\subseteq\mathbb{R}$ be its value group. We investigate the power series ring $V[[t]]$, and finally answer the long standing problem of when this ring is coherent. Specifically, we prove that it is coherent if and only if $V$ is \emph{almost maximal}, and either $\Gamma=\mathbb{R}$ or $\Gamma$ is discrete. Furthermore, we demonstrate that the weak Hilbert--syzygy identity $\wdim V[[X]]=\wdim V+1$ holds if and only if $V[[X]]$ is a coherent ring.
\end{abstract}

\maketitle

\section{Introduction}

\noindent In this article, we consider all rings to be commutative with identity. Recall that a ring $S$ is \emph{coherent} if every finitely presented $S$-module has a finite free resolution (or equivalently, every finitely generated ideal of $S$ is finitely presented). 

One of the most classical examples of a coherent ring is a \emph{valuation domain}, i.e. commutative domain $V$ where for all ideals $I,J\trianglelefteq V$, either $I\subseteq J$ or $J\subseteq I$. In this case, all finitely generated ideals are principal. Note that all fractional ideals of $K=Frac(V)$ can be indexed by a totally ordered abelian group $\Gamma$, which we call the \emph{value group} of $V$.

It is clear that any Noetherian ring is coherent. But unlike Noetherianity, it it not generally true that constructions such as polynomial rings or power series rings preserve coherence; there are examples of coherent rings $S$ such that $S[X]$ is not coherent.

However, if $V$ is a valuation ring, then $V[X]$ is a coherent ring, but the power series ring $V[[X]]$ is not in general. In fact, it was proved by J{\o}ndrup and Small in \cite[Section 1]{JondrupSmall} that if $V$ has rank greater than 1, then $V[[X]]$ is never coherent, and they remarked that it is unknown whether or not $V[[X]]$ is coherent when $V$ has rank 1 (i.e. when $\Gamma$ is a subgroup of $\mathbb{R}$).

All subgroups of $\mathbb{R}$ are either discrete or dense. If $\Gamma$ is discrete, then $V$ is Noetherian, so $V[[X]]$ is also Noetherian. However, it was first proved in \cite{coherence1}, and later by the second author in \cite{Jones} that if $\Gamma$ is proper and dense in $\mathbb{R}$, then $V[[X]]$ is not coherent, and it was recently proved by the third author in \cite{Wang} that if $\Gamma=\mathbb{R}$ and $V$ is complete, then $V[[X]]$ is coherent if and only if $K$ is spherically complete. In general, when $\Gamma=\mathbb{R}$, the question of coherence of $V[[X]]$ remains open.

In this article, we link this question to another open question for power series rings, namely the question of Hilbert's syzygy theorem:

For a ring $S$ and an $S$-module $M$, let $\fd_S M$ denote the flat dimension of $M$, and let $\wdim S$ denote the weak global dimension of $S$, i.e $$\wdim S:=\sup\{\fd_S(M):M\text{ an }S{-module}\}$$

The classical Hilbert's syzygy theorem states that for any ring $R$, $\wdim R[X]=\wdim R+1$, and an ongoing project in commutative algebra is to establish for which rings $R$ the same result holds for the power series ring, i.e.
\begin{equation}\label{eq:coherent-equality}
  \wdim R[[X]]=\wdim R+1
\end{equation}

\noindent For example, $V$ is a non-field valuation ring, then all ideals in $V$ are flat, so $\wdim(V)=1$. Thus $V[[X]]$ satisfies Hilbert's syzygy identity \eqref{eq:coherent-equality} if and only if $\wdim(V[[X]])=2$.\\

\noindent It was proved by J{\o}ndrup and Small that if the power series ring $R[[X]]$ is coherent (e.g. if $R$ is Noetherian), then \eqref{eq:coherent-equality} is satisfied \cite[Lemma 1]{JondrupSmall}. The question was investigated further by Bouchiba in \cite{Bouchiba}, who claimed the identity holds for any coherent ring $R$, and moreover that the stronger identity $\fd_{R[[X]]}M=\fd_R(M)+1$ holds for any $R[[X]]$-module $M$. 

However, the second author identified a gap in this argument in \cite{Bouchiba}, and in fact showed that if $V$ is a rank-one valuation domain $V$ with proper dense value group \cite[Theorem A]{Jones}. More precisely, there exists a $V[[X]]$-module $C$ satisfying
\[
  \fd_V C=0,
  \qquad
  \fd_{V[[X]]}C\geq2.
\]
It was remarked in the conclusion to \cite{Jones} that this does not itself disprove Hilbert's syzygy identity of $V$, and that to this end an additional homological step would be required. Thus the question remained open.

The main result of this article provides a full answer to both these open questions. In the below, we let $V((X)):=V[[X]][X^{-1}]$ be the ring of Laurent series over $V$, and we define what it means for $V$ to be \emph{almost maximal} in section 2.



\begin{letterthm}\label{main_theorem}
Let $V$ be a valuation ring of rank 1, with value group $\Gamma\subseteq\mathbb{R}$. Then the following are equivalent:

\begin{itemize}
    \item[i.] $V[[X]]$ is a coherent ring.

    \item[ii.] $V[[X]]$ satisfies Hilbert's syzygy identity \eqref{eq:coherent-equality}.

    \item[iii.] $V$ is an almost maximal valuation ring, and either $\Gamma=0$, $\Gamma=\mathbb{R}$ or $\Gamma\cong\mathbb{Z}$.
\end{itemize}

\noindent In all other cases, $\wdim V[[X]]=\wdim V((X))>2$.

\end{letterthm}

\noindent\textbf{Note:} The condition $\Gamma=0$ just means that $V$ is a field, while $\Gamma\cong\mathbb{Z}$ means that $V$ is discrete. In both these cases, $V$ is Noetherian, so we know that $V[[X]]$ is Noetherian and satisfies $\wdim(V[[X]])=\gdim(V[[X]])=\gdim(V)+1$.\\

\noindent Again, $i \implies ii$ was proved by  J{\o}ndrup-Small in \cite[Lemma 1]{JondrupSmall}, so this theorem constitutes a partial converse to their result. To complete the proof, the only steps remaining are to prove $i\iff iii$ and $ii\implies iii$.\\

\noindent After proving some essential preliminary results in section 2, we go on in section 3 to complete the proof that $V[[X]]$ is a coherent ring if and only if $V$ satisfies the stated conditions in $iii$. The proof that $i\implies iii$ was given by the third author in \cite[Section 5]{Wang} (with the assumption of maximality rather than almost maximality, but the argument still holds by \cite[Remark 4.3]{Wang}), and the proof of $iii\implies i$ is completely analogous to the argument in \cite[Section 3]{Wang}. \\

\noindent The remainder of the paper is dedicated to proving $ii\implies iii$. Specifically, we will prove that if $V$ does not satisfy the conditions stated in $iii$, then $V[[X]]$ does not satisfy Hilbert's syzygy identity. Moreover, the inequality $\wdim R[[X]]\geq\wdim R+1$ is satisfied for all coherent rings $R$ by the proof of \cite[Lemma 1.1]{JondrupSmall}, so since $\wdim V=1$ for any non-field valuation ring $V$, we only need to prove that $\wdim V[[X]]\geq 3$ if $V$ does not satisfy the stated conditions.


To show that $\wdim R[[X]]\geq 3$ for any commutative domain $R$, it suffices to show that there exist elements $f,r\in R[[X]]$ such that $(f)\cap (r)$ is not flat in $R[[X]]$, and apply two dimensional Tor-shifts. So it remains to construct such elements $f,r\in V[[X]]$ under the appropriate conditions on $V$.

In section 4, we consider the case where $\Gamma$ is a proper, dense subgroup of $\mathbb{R}$, and following the approach of the second author in \cite{Jones} we utilise the coefficient-infimum function
\begin{equation}\label{eq:intro-valuation}
  \nu\!\left(\sum_{n\geq0}a_nX^n\right)
  =\inf_{n\geq0}v(a_n)
\end{equation}
which is a valuation on $V[[X]]$ by \cite[Corollary~2.5]{Jones}. We consider the cut ideals
\[
  I_s=\{g\in V[[X]]:\nu(g)\geq s\}.
\]
for each $s>0$. A fractional-intersection identity satisfied by every nonzero flat ideal then shows that $I_c$ is not flat if $c\notin\Gamma$. Thus we construct $f,r\in V[[X]]$ such that $(f)\cap(r)=fI_c$, and hence $(f)\cap (r)$ is not flat as required.

In section 5, we turn to the case where $\Gamma=\mathbb{R}$, and focus on the case where $K=Frac(V)$ is complete, but not spherically complete. This was the case considered by the third author in \cite[Section 4]{Wang}, where a pair $f,r\in V[[X]]$ was constructed from the tails of infinite products, and the intersection $(f)\cap (r)$ was proved to be of infinite presentation. We strengthen this construction homologically, and apply a flat-descent lemma to show that, under a suitable coprimality condition in a GCD overring, flatness of the intersection would force invertibility. This suffices to prove that $(f)\cap (r)$ is not flat.

Finally, in section 6, we consider the remaining case where $\Gamma=\mathbb{R}$ but $V$ is not almost maximal. Equivalently, the topological completion $\widehat{K}$ is not spherically complete. Setting $\widehat{V}=\mathcal{O}_{\widehat{K}}$, we prove that $\wdim V[[X]] = \wdim\widehat{V}[[X]]$, and hence $\wdim(V[[X]])>2$.\\

\noindent The question of whether Hilbert's syzygy identity \eqref{eq:coherent-equality} is satisfied for any valuation ring $V$ of rank greater than 1 remains open. However, given that no example of any $R$ such that $R[[X]]$ is not coherent but the syzygy identity is satisfied seems to appear in the literature, we tentatively make the following conjecture.

\begin{conjecture}
If $R$ is a coherent ring, then $R[[X]]$ is coherent if and only if $\wdim(R[[X]])=\wdim(R)+1$.
\end{conjecture}

\section{Prerequisites}

\noindent In this section, we will prove some technical results which will be essential in our main argument.

\subsection{Flatness criteria}

For now, let $A$ be any commutative domain (later we will realise $A$ as a power series ring). We now isolate some key flatness critera.

\begin{lemma}\label{lem:flat-intersection}
Let $L=Frac(A)$, and let $E\neq 0$ be a flat
ideal of $A$. If $P,Q\subseteq L$ are $A$-submodules, then
\begin{equation}\label{eq:flat-intersection}
  E(P\cap Q)=EP\cap EQ.
\end{equation}
\end{lemma}

\begin{proof}
Tensor the exact sequence
\[
  0\longrightarrow P\cap Q
  \longrightarrow P\oplus Q
  \longrightarrow P+Q
  \longrightarrow0
\]
with $E$. For every $A$-submodule $P\subseteq L$, the multiplication
map
\[
  E\otimes_AP\longrightarrow L
\]
is injective. Indeed, flatness preserves the injection $P\hookrightarrow
L$, and, with $U=A\setminus\{0\}$,
\[
  E\otimes_AL\cong U^{-1}E=L
\]
because $E\neq 0$. Under the multiplication map, the image of
$E\otimes_AP$ is $EP$. The kernel in the tensored sequence therefore
identifies both with $E(P\cap Q)$ and with $EP\cap EQ$, proving
\eqref{eq:flat-intersection}.
\end{proof}

\begin{lemma}\label{thm:main}
Suppose $f,r\in A$ and the ideal $(f)\cap (r)$ is not a flat $A$-module. Then
\[
  \fd_A\bigl(A/(f,r)\bigr)\geq3.
\]
Consequently,
\begin{equation}\label{eq:strict-global}
  \wdim A\geq 3
\end{equation}
\end{lemma}

\begin{proof}
Let
\[
  L=(f)\cap(r),
  \qquad
  J=(f,r).
\]
By assumption, $L$ is not flat. The
standard exact sequence
\begin{equation}\label{eq:first-shift}
  0\longrightarrow L
  \longrightarrow(f)\oplus(r)
  \longrightarrow J
  \longrightarrow0
\end{equation}
has flat middle term, since $(f)\cong A\cong(r)$. Hence there is an
$A$-module $N$ such that $\Tor_1^A(L,N)\neq0$, and the long exact
$\Tor$ sequence associated with \eqref{eq:first-shift} gives
\[
  \Tor_2^A(J,N)\cong\Tor_1^A(L,N)\neq0.
\]
A second dimension shift, applied to
\[
  0\longrightarrow J\longrightarrow A
  \longrightarrow A/J\longrightarrow0,
\]
yields
\[
  \Tor_3^A(A/J,N)\cong\Tor_2^A(J,N)\neq0.
\]
Therefore $\fd_A(A/J)\geq3$.\end{proof}

\begin{lemma}[Flat descent criterion]\label{lem:flat-descent}
Let $B$ be a domain with $A\subseteq B$, viewed inside $\Frac(B)$, and let $E$ be a non-zero ideal of $A$. We also assume
\begin{itemize}
\item[(i)] $B$ is a GCD domain,

\item[(ii)] $\Frac(A)\cap B=A$,

\item[(iii)] There exist $s,y,z\in B$ such that $E\subseteq sB$, $sy,sz\in E$ and $\underset{B}{\text{\emph{gcd}}}(y,z)=1$.
\end{itemize}

\noindent Then if $E$ is flat over $A$, it follows that $E$ is invertible.  In particular, $E$ is finitely generated.
\end{lemma}

\begin{proof}
Put $u=sy$ and $w=sz$.  Then $u,w\in E\subseteq A$.  For a nonzero
ideal $I$ of $A$, write
\[
  I^{-1}=(A:I)
  =\{q\in\Frac(A):qI\subseteq A\}.
\]
We claim that
\begin{equation}\label{eq:inverse-two-generators}
  E^{-1}=u^{-1}A\cap w^{-1}A.
\end{equation}
The inclusion from left to right follows from $u,w\in E$.  Conversely,
let $q$ belong to the right-hand side and set $x=qs\in\Frac(B)$.  Then
\[
  xy=qu\in A\subseteq B,
  \qquad
  xz=qw\in A\subseteq B.
\]
For relatively prime elements of a GCD domain one has
$yB\cap zB=yzB$.  Hence
\[
  xyz\in yB\cap zB=yzB,
\]
and cancellation gives $x=qs\in B$.  If $e\in E$, write $e=sb$ with
$b\in B$.  Then
\[
  qe=(qs)b\in B\cap\Frac(A)=A
\]
Thus $q\in E^{-1}$, proving \eqref{eq:inverse-two-generators}.

Assume that $E$ is flat.  Lemma~\ref{lem:flat-intersection}, applied to
$P=u^{-1}A$ and $Q=w^{-1}A$, gives
\[
  EE^{-1}=E(u^{-1}A)\cap E(w^{-1}A).
\]
Both ideals on the right contain $1$, since
$1=u\cdot u^{-1}=w\cdot w^{-1}$.  Therefore $EE^{-1}=A$, so $E$ is
invertible.  If
\[
  1=\sum_{j=1}^m e_jq_j,
  \qquad e_j\in E,\quad q_j\in E^{-1},
\]
then every $e\in E$ satisfies
$e=\sum_j e_j(q_je)$ with $q_je\in A$.  Hence
$E=(e_1,\ldots,e_m)$.
\end{proof}


\noindent We also record the following known consequence of a regular-element theorem of Rajabi, Torrecillas and Yassemi, and which will complete the proof of the final statement of our main theorem. 

\begin{proposition}\label{prop:known-localization}
For every commutative ring $R$,
\begin{equation}\label{eq:known-localization}
  \wdim R[[X]]
  =\max\{\wdim R+1,\wdim R((X))\},
\end{equation}
where $R((X))=R[[X]][X^{-1}]$.
\end{proposition}

\begin{proof}
Put $S=R[[X]]$. The element $X$ is a central nonzerodivisor of $S$,
$S/XS\cong R$, and $S_X=R((X))$. If $\wdim R<\infty$, the assertion is
exactly \cite[Theorem~3.7]{RajabiTorrecillasYassemi}. If
$\wdim R=\infty$, the unconditional lower bound
$\wdim R[[X]]\geq\wdim R+1$ from the proof of
\cite[Lemma~1]{JondrupSmall} gives $\wdim R[[X]]=\infty$, and
\eqref{eq:known-localization} follows as well.
\end{proof}

\subsection{Almost maximal valuation rings}

We will now fix a valuation ring $V$ of rank 1, with fraction field $K$ and associated valuation $v:K^\times\to \mathbb{R}$. We write $\mathfrak{m}$ for the maximal ideal of $V$.\\

\noindent Recall that a valued field $K$ is \emph{spherically complete} if any descending chain of closed balls in $K$ has non-empty intersection (which implies that $K$ is complete). If $K=Frac(V)$ is spherically complete then we say $V$ is a \emph{maximal} valuation ring.

Furthermore, we say that $K$ is \emph{almost spherically complete} if for any $R>0$ and any descending chain of closed balls in $K$ of radius at least $R$, the intersection of this chain is non-empty. Any spherically complete field $K$ is almost spherically complete, and if $K=Frac(V)$ is almost spherically comlete, then we say $V$ is an \emph{almost maximal valuation ring}.\\

\noindent\emph{Example:} Any discrete valuation ring is almost maximal, and the ring of Hahn series over any totally ordered abelian group is maximal.\\

\noindent The following elementary result will be essential when alternating between almost maximal valuation rings and their completions.

\begin{lemma}\label{lem: maximal completion}
For any valuation ring $V$ of rank 1, the completion $\widehat{V}$ of $V$ is the valuation ring of $\widehat{K}$. Moreover, $V$ is almost maximal if and only if $\widehat{V}$ is maximal.
\end{lemma}

\begin{proof}

Since the valuation $v$ extends to $\widehat{K}$, it is clear that $V\subseteq\mathcal{O}_{\widehat{K}}=\{x\in\widehat{K}:v(x)\geq 0\}$. Moreover, suppose $x\in\mathcal{O}_{\widehat{K}}$, then $x$ is the limit of a sequence $(x_n)$ with each $x_n\in V$, and hence for sufficiently $n$, $v(x_n)=v(x)\geq 0$, thus $x\in\widehat{V}$ as required.\\

\noindent For the second statement, note that a closed ball in $K$ has the form $y+I_\alpha$ or $y+I_{\alpha}^+$, where $I_{\alpha}=\{x\in K:v(x)\geq\alpha\}$ and $I_{\alpha}^+=\{x\in K:v(x)>\alpha\}$. We can take the radius of the ball to be $2^{-\alpha}$. Closed balls in $\widehat{K}$ have a similar form.

But for any $y\in\widehat{K}$, by density there exists $y'\in K$ such that $v(y-y')>\alpha$, so $y+I_{\alpha}=y'+I_{\alpha}$ and $y+I_{\alpha}^+=y'+I_{\alpha}^+$, so a closed ball in $\widehat{K}$ contains a closed ball in $K$ with the same centre and radius.\\

\noindent Therefore, a descending sequence of closed balls with bounded radius has non-empty intersection in $K$ if and only if it has non-empty intersection in $\widehat{K}$, i.e. $K$ is almost spherically closed if and only if $\widehat{K}$ is. 

Moreover, a descending sequence of closed balls whose radii converge to 0 corresponds to a Cauchy sequence, so naturally it converges in $\widehat{K}$ by completeness. Thus $\widehat{K}$ is spherically complete if and only if $K$ is almost spherically complete as required.\end{proof}

\noindent Almost maximal valuation rings were studied by Matlis in \cite{MatlisPruefer}. There, it was defined to be a valuation ring $V$ such that $V/I$ is linearly compact for all nonzero ideals $I$. Here, a module is \emph{linearly compact} if every family of cosets of submodules with the finite intersection property has nonempty intersection. It is well known that finite direct sums and quotients of linearly compact modules are
linearly compact.

It is immediately verified that these definitions are equivalent, but the formulation in terms of linearly compact modules will be essential in the proof of the following technical results (analagous to \cite[Lemma 2.1 and Lemma 2.2]{Wang}), and in the proof of our main theorem.


\begin{lemma}\label{am:affine}
Let $F=V^n$ be a finite free $V$-module, and let
$Z_1\supseteq Z_2\supseteq\cdots$ be nonempty cosets of submodules $L_i$
of $F$. If $cF\subseteq L_i$ for every $i$, for some $0\ne c\in V$,
then $\bigcap_i Z_i\ne\varnothing$.
\end{lemma}

\begin{proof}
The images of the $Z_i$ in the linearly compact module $F/cF\cong (V/cV)^n$ have
nonempty intersection. Every lift of a point in this intersection
belongs to all the $Z_i$, since $cF\subseteq L_i$.
\end{proof}

\begin{proposition}\label{am:lattices}
If $v(K^\times)=\mathbb{R}$, then the intersection of any nonempty family of finitely generated $V$-submodules of a finite free $V$-module is finite free.
\end{proposition}

\begin{proof}
We adapt the rank induction of
\cite[Lemma~6.14]{BhattScholzeWittGrassmannian}, using Matlis's theorem
\cite[Theorem~4]{MatlisPruefer}, which states that the module $K/V$ is injective.

The module $K$ is also injective: a homomorphism from a nonzero ideal
of $V$ to $K$ is multiplication by an element of $K$, and hence
extends to $V$. Thus applying the functor Hom$_V(V/\mathfrak{m},-)$ to the sequence $0\to V\to K\to K/V\to0$ and taking the associated long-Ext exact sequence we see that $$\operatorname{Ext}^2_V(V/\mathfrak m,V)=0$$  and applying dimension
shifting gives
\[
\operatorname{Ext}^1_V(\mathfrak m,V)
\simeq \operatorname{Ext}^2_V(V/\mathfrak m,V)=0.
\]

\noindent For a submodule $N\subseteq K^d$, identify
$\operatorname{Hom}_V(\mathfrak m,N)$ with
\[
N^\sharp=\{z\in K^d:\mathfrak m z\subseteq N\}.
\]
Indeed, a homomorphism extends upon localization to a $K$-linear map
$K\to K^d$, so is multiplication by a unique vector.
We claim that a bounded submodule $N\subseteq K^d$ satisfying
$N=N^\sharp$ is finite free. Here bounded means $N\subseteq aV^d$
for some $a\in K^\times$.

For $d=1$, a nonzero bounded submodule has a real lower bound
$\gamma=\inf v(N\setminus\{0\})$ and is either
$\{z:v(z)\ge\gamma\}$ or $\{z:v(z)>\gamma\}$, with zero included.
The second possibility contradicts $N=N^\sharp$; the first is
principal because $v(K^\times)=\mathbb R$. The zero module causes
no difficulty.

For the induction step, let $N'=N\cap K^{d-1}$ and let $N''\subseteq K$
be the image of $N$ under the last-coordinate projection. Both are
bounded, and $(N')^\sharp=N'$. Induction makes $N'$ finite free, so
$\operatorname{Ext}^1_V(\mathfrak m,N')=0$. Apply
$\operatorname{Hom}_V(\mathfrak m,-)$ to
$0\to N'\to N\to N''\to0$.
The natural multiplication maps identify the first two modules with
their Hom modules; exactness then gives $(N'')^\sharp=N''$.
The rank-one case makes $N''$ finite free. The sequence splits, proving
the claim.

Now let $L=\bigcap_i L_i\subseteq V^d$, with each $L_i$ finitely generated. it follows from \cite[Lemma 3.1.1]{AntieauDatta} that $L_i$ is finite free, so since the value group is $\mathbb R$, one has $V^\sharp=V$:
if $v(z)<0$, choose $t\in\mathfrak m$ with $0<v(t)<-v(z)$.
Using a basis of $L_i$ in its $K$-span gives $L_i^\sharp=L_i$.
Consequently
\[
L^\sharp=\bigcap_i L_i^\sharp=\bigcap_i L_i=L.
\]
The claim applies to the bounded module $L$.
\end{proof}

\section{Coherence over an almost maximal valuation domain}
\label{am:section}

\noindent In this section, we will answer J{\o}ndrup and Smalls question from \cite{JondrupSmall}, and once and for all determine when a power sereis ring over a valuation ring is coherent. Again, fix an almost maximal valuation ring $V$ of rank 1, with associated valuation $v:K^\times\to\mathbb{R}$, and we will assume that $v(K^\times)=\mathbb{R}$.

Let $A:=V[[X]]$. We will adapt the kernel-lifting argument of \cite[Section~3]{Wang} to prove that $A$ is a coherent ring.\\


\noindent For each $r\in\mathbb{N}$, let $A_r=A/X^rA$. As a $V$-module, $A_r\cong V^r$ is finite free.
For a finite matrix $\Phi:A^m\to A^n$, write
\[
\mathcal K=\ker\Phi,\qquad
\mathcal K_r=\ker(\Phi_r:A_r^m\to A_r^n),
\]
and let $\rho_{s,r}:\mathcal K_s\to\mathcal K_r$ denote reduction,
for $s\ge r$. Each $\mathcal K_r$ is finite free over $V$.  Completeness of $A$ gives
\begin{equation}\label{am:limit}
\mathcal K\simeq\varprojlim_r\mathcal K_r.
\end{equation}
Set
\[
E_r=\bigcap_{s\ge r}\operatorname{im}\rho_{s,r}.
\]
An element of $E_r$ lifts to every finite level separately; the next
proposition supplies compatible lifts.

The ideal $(X)$ of $A$ is prime. Every nonzero element of
$D=A_{(X)}$ is a power of $X$ times a unit, so $D$ is a DVR with
uniformizer $X$ and residue field $K$. Moreover
\[
D_r=D/X^rD\simeq K[X]/(X^r)\simeq A_r\otimes_V K.
\]
Localization is flat, and hence
\begin{equation}\label{am:generic}
\mathcal K_D:=\ker(\Phi_D:D^m\to D^n)=\mathcal K\otimes_A D.
\end{equation}
As a submodule of a finite free module over a DVR, $\mathcal K_D$
is finite free. It can therefore be generated over $D$ by finitely
many vectors coming from $\mathcal K$; this does not presuppose
finite generation of $\mathcal K$.

\begin{proposition}\label{am:lifting}
Every map $E_{r+1}\to E_r$ is surjective. In particular,
\[
\operatorname{im}(\mathcal K\to\mathcal K_1)
=\bigcap_{r\ge1}\operatorname{im}(\mathcal K_r\to\mathcal K_1).
\]
\end{proposition}

\begin{proof}
Fix $x_r\in E_r$. For $s\ge r+1$, put
\[
Y_s=\operatorname{im}\rho_{s,r+1},\qquad
Z_s=\{z\in Y_s:\rho_{r+1,r}(z)=x_r\},
\qquad N_r=\ker\rho_{r+1,r}.
\]
The $Z_s$ form a descending sequence of nonempty cosets with directions
$L_s=Y_s\cap N_r$ (i.e. $Z_s=z+L_s$ for any $z\in Z_s$). The modules $Y_s$ and $N_r$ are finitely generated
$V$-submodules of $A_{r+1}^m$, so their intersections $L_s$ are finite
free by coherence of $V$. Also $N_r$ lies in the coefficient layer
$X^rA^m/X^{r+1}A^m\simeq V^m$.

Tensoring by $K$ is exact, so it preserves these kernels, images and
intersections. Choose Smith normal form for $\Phi_D$, with nonzero
diagonal entries $X^{a_1},\ldots,X^{a_k}$. Let
$a=\max_i a_i$ when $k>0$, and $a=0$ otherwise. In the corresponding
source basis $e_1,\ldots,e_m$, for $s\ge\max\{1,a\}$,
\[
\ker\Phi_{D_s}
=\bigoplus_{i=1}^k X^{s-a_i}D_s e_i
 \ \oplus\ \bigoplus_{i=k+1}^m D_s e_i.
\]
Thus, for $s\ge r+1+a$, the image of this kernel modulo $X^{r+1}$
equals the image of $\mathcal K_D$. Write
\[
q:D^m\to D_{r+1}^m,\qquad
H=X^rD^m/X^{r+1}D^m.
\]
For these $s$, flat base change consequently gives
\[
L_s\otimes_V K=q(\mathcal K_D)\cap H.
\]
Since $q^{-1}(H)=X^rD^m$ and $\mathcal K_D$ is $X$-saturated,
\begin{equation}\label{am:directions}
q(\mathcal K_D)\cap H
=q(\mathcal K_D\cap X^rD^m)
=q(X^r\mathcal K_D)=:W.
\end{equation}

Choose $u_1,\ldots,u_t\in\mathcal K$ generating $\mathcal K_D$ over
$D$, using \eqref{am:generic}. Let $L_0\subseteq N_r$ be the $V$-span
of the classes of $X^r u_j$. These are reductions of actual kernel
vectors at every level, and \eqref{am:directions} gives
\[
L_0\subseteq L_s\quad\text{for all }s,\qquad
L_0\otimes_V K=L_s\otimes_V K=W\quad(s\ge r+1+a).
\]
Choose $s_0\ge r+1+a$. The finite modules $L_0$ and $L_{s_0}$
have the same $K$-span. Clearing denominators for generators of
$L_{s_0}$ gives $0\ne c\in V$ such that
\[
cL_{s_0}\subseteq L_0\subseteq L_s\subseteq L_{s_0}
\qquad(s\ge s_0).
\]
Choose $z_0\in Z_{s_0}$. The translated sets $Z_s-z_0$, for $s\ge s_0$,
are nested nonempty cosets in the finite free module $L_{s_0}$, with
directions containing $cL_{s_0}$. Applying Lemma~\ref{am:affine}, there exists a point in their intersection. Translating back gives an element of
$E_{r+1}$ mapping to $x_r$.

Starting from any $x_1\in E_1$, recursively choose compatible
$x_r\in E_r$. Equation~\eqref{am:limit} produces the required global
kernel vector. The reverse inclusion in the asserted equality follows
by reduction of a global vector.
\end{proof}

\begin{theorem}\label{thm: coherence}
$A$ is a coherent domain.    
\end{theorem}

\begin{proof} For each $r$, the module
$\Lambda_r=\operatorname{im}(\mathcal K_r\to\mathcal K_1)$
is a finitely generated $V$-submodule of $V^m$.
Proposition~\ref{am:lifting} and Proposition~\ref{am:lattices} show that
\[
\operatorname{im}(\mathcal K\to\mathcal K_1)=\bigcap_r\Lambda_r
\]
is finite free. The kernel $\mathcal K$ is $X$-saturated in $A^m$,
so
\[
\mathcal K/X\mathcal K
\simeq\operatorname{im}(\mathcal K\to\mathcal K_1).
\]
It is closed in the complete module $A^m$, and saturation gives
$\mathcal K\cap X^rA^m=X^r\mathcal K$. Hence $\mathcal K$ is complete
and separated for its $X$-adic topology. Lifts of finitely many
generators of $\mathcal K/X\mathcal K$ generate $\mathcal K$:
successive subtraction of their linear combinations and division by
$X$ expresses any vector using convergent power-series coefficients.
Thus every finite matrix over $A$ has finitely generated kernel, and
$A$ is coherent.\end{proof}

\begin{corollary}\label{cor: conditions for coherence}
Let $V$ be any valuation ring. Then the power series ring $V[[X]]$ is coherent if and only if

\begin{itemize}
    \item $V$ has rank 1,

    \item $V$ is almost maximal, and

    \item $V$ is discrete, or has value group $\mathbb{R}$.
\end{itemize}
\end{corollary}

\begin{proof}
Firstly, if $V$ has rank 1, then $V[[X]]$ is Noetherian, and hence coherent, if $V$ is discrete. On the other hand, if $V$ is almost maximal with value group $\mathbb{R}$, then $V[[X]]$ is coherent by Theorem \ref{thm: coherence}.\\

\noindent For the converse, if $V$ has rank greater than 1, then we know that $V[[X]]$ is not coherent by \cite[Section 1]{JondrupSmall}. So we may assume that $V$ has rank 1, and let $v:K^\times\to\mathbb{R}$ be the associated valuation.

Since the value group $\Gamma=v(K^\times)$ is a subgroup of $\mathbb{R}$, either $\Gamma$ is discrete, or $\Gamma$ is dense in $\mathbb{R}$. If $\Gamma\neq\mathbb{R}$ and $\Gamma$ is dense, it follows from \cite[Section 3 Corollary]{coherence1} and \cite[Theorem 3.1]{Jones} that $V[[X]]$ is not coherent. 

Finally, if $\Gamma=\mathbb{R}$ but $V$ is not almost maximal, then the proof in \cite[Section 4 and 5]{Wang} shows that $V[[X]]$ is not coherent.\\

\noindent\textbf{Note:} The proof in \cite{Wang} assumed that $V$ was complete, but the only place in section 4 or 5 that uses completeness is the proof of \cite[Construction Lemma 4.2]{Wang}, where it was required to assert that a descending chain of balls with empty intersection has bounded radius (as noted in \cite[Remark 4.3]{Wang}). Since assuming $V$ is not almost maximal implies the existence of such a chain immediately, the proof goes through without any changes.\end{proof}

\section{Cut ideals and flatness}

\noindent We now turn our attention to Hilbert's syzygy. Again, fix a valuation domain $V$ of rank 1, with fraction field $K$, valuation $v:K^\times\to\mathbb{R}$, and let $\Gamma=v(K^\times)\subseteq\mathbb{R}$ be the value group.

In this section, we will assume that $\Gamma$ is a proper, dense subgroup of $\mathbb{R}$. Set $A:=V[[X]]$, and consider the valuation $\nu$ on $A$ given by \eqref{eq:intro-valuation} (and also denote by $\nu$ its extension to $\Frac(A)$).\\

\noindent For $s\geq0$, set
\[
  I_s=\{a\in A:\nu(a)\geq s\}.
\]

Since $\nu(a)\geq0$ for $a\in A$, each $I_s$ is an ideal of $A$. Moreover, if $s\in \Gamma$ then $s=v(a)$ for some $a\in V$, so $I_s=(a)$ is a principal ideal, and hence flat. Fixing $c>0$ with $c\notin\Gamma$, we will first prove that $I_c$ is not a flat ideal of $A$.

Choose $\rho\in\Gamma$ with $\rho>c$, and set
\[
  d=\rho-c.
\]
Then
\begin{equation}\label{eq:cuts}
  c,d>0,
  \qquad c,d\notin\Gamma,
  \qquad c+d=\rho.
\end{equation}
Choose $r\in V$ with $v(r)=\rho$. By density, there are sequences
$(c_n)$ and $(d_n)$ in $\Gamma$ such that
\[
  c_n>c,\quad c_n\longrightarrow c,
  \qquad
  d_n>d,\quad d_n\longrightarrow d.
\]
Choose $u_n,w_n\in V$ with $v(u_n)=c_n$ and $v(w_n)=d_n$, and set
\begin{equation}\label{eq:hf}
  h=\sum_{n\geq0}u_nX^n,
  \qquad
  f=\sum_{n\geq0}w_nX^n.
\end{equation}
Then
\begin{equation}\label{eq:hf-values}
  \nu(h)=c,
  \qquad
  \nu(f)=d.
\end{equation}
Moreover,
\begin{equation}\label{eq:Irho}
  I_\rho=rA.
\end{equation}
Indeed, $\nu(g)\geq\rho=v(r)$ if and only if every coefficient of $g$
is divisible by $r$ in $V$.

\begin{lemma}\label{lem:strict-product}
With the notation above,
\[
  I_cI_d\subsetneq I_\rho.
\]
\end{lemma}

\begin{proof}
Multiplicativity of $\nu$ gives $I_cI_d\subseteq I_\rho$. If
$a=\sum a_nX^n\in I_c$, every nonzero coefficient $a_n$ has value
strictly greater than $c$, because $v(a_n)\in\Gamma$ while
$c\notin\Gamma$. The analogous assertion holds for elements of $I_d$. Hence the constant coefficient of every finite sum of products of an
element of $I_c$ and an element of $I_d$ is either zero or has value
strictly greater than $c+d=\rho$. The constant series $r$ lies in
$I_\rho$ and has constant coefficient of value exactly $\rho$.
Therefore $r\notin I_cI_d$.
\end{proof}

\begin{proposition}\label{prop:cut-not-flat}
The ideal $I_c$ is not flat over $A$.
\end{proposition}

\begin{proof}
Put
\[
  q=\frac{h}{r}\in\Frac(A).
\]
By \eqref{eq:cuts} and \eqref{eq:hf-values},
\[
  \nu(q)=c-\rho=-d.
\]
We first compute
\begin{equation}\label{eq:A-qA}
  A\cap qA=qI_d.
\end{equation}
If $qa\in A$, then
\[
  0\leq\nu(qa)=\nu(a)-d,
\]
so $a\in I_d$. Conversely, if $a\in I_d$, then
$\nu(ha)\geq c+d=\rho$. Every coefficient of $ha$ is divisible by
$r$, and therefore $qa=ha/r\in A$. This proves
\eqref{eq:A-qA}.

We also have
\begin{equation}\label{eq:Ic-qIc}
  I_c\cap qI_c=qI_\rho.
\end{equation}
Indeed, if $z=qa\in I_c$ with $a\in I_c$, then
\[
  \nu(a)=\nu(z)+d\geq c+d=\rho,
\]
so $a\in I_\rho$. Conversely, if $a\in I_\rho$, then
$a\in I_d$, whence $qa\in A$ by \eqref{eq:A-qA}; moreover,
$a\in I_c$ and $\nu(qa)\geq\rho-d=c$. Thus
$qa\in I_c\cap qI_c$, proving \eqref{eq:Ic-qIc}.

Suppose that $I_c$ were flat. Lemma~\ref{lem:flat-intersection},
applied to $E=I_c$, $P=A$ and $Q=qA$, would give
\[
  I_c(A\cap qA)=I_c\cap qI_c.
\]
Using \eqref{eq:A-qA} and \eqref{eq:Ic-qIc}, we obtain
\[
  qI_cI_d=qI_\rho,
\]
and hence $I_cI_d=I_\rho$, contrary to
Lemma~\ref{lem:strict-product}.
\end{proof}

\begin{lemma}\label{lem:principal-intersection}
For the elements $f$ and $r$ chosen above,
\begin{equation}\label{eq:principal-intersection}
  (f)\cap(r)=fI_c.
\end{equation}
In particular, $(f)\cap(r)$ is not flat over $A$.
\end{lemma}

\begin{proof}
Let $fg=ra$ with $g,a\in A$. By multiplicativity of $\nu$,
\[
  d+\nu(g)=\rho+\nu(a)\geq\rho,
\]
so $\nu(g)\geq\rho-d=c$ and $g\in I_c$. Thus
$(f)\cap(r)\subseteq fI_c$.

Conversely, if $g\in I_c$, then
\[
  \nu(fg)\geq d+c=\rho.
\]
Every coefficient of $fg$ is divisible by $r$, and hence
$fg\in(r)$. This proves \eqref{eq:principal-intersection}. Since
$A$ is a domain and $f\neq0$, multiplication by $f$ identifies $I_c$
with $fI_c$. Proposition~\ref{prop:cut-not-flat} gives the final
assertion.
\end{proof}

\begin{theorem}\label{thm:main'}
If $\Gamma=v(K^\times)$ is a proper dense subgroup of $\mathbb{R}$, then $\wdim V[[X]]>2$.
\end{theorem}

\begin{proof}

Using Lemma \ref{lem:principal-intersection} we know that there exist $f,r\in V[[X]]$ such that $(f)\cap (r)$ is not flat as a $V[[X]]$-module, and hence $\wdim V[[X]]\geq 3$ by Theorem \ref{thm:main}.\qedhere 

\end{proof}

\section{The full real value group}\label{sec:full-real}

\noindent We now return to valued fields with full real value group. The cut argument of Sections~2--3 is unavailable in this case, since there is no real cut outside the value group.

Let $K$ be a complete valued field with
$v(K^\times)=\mathbb R$, let $V=\mathcal O_K$, and put $A=V[[X]]$.
Assume that $K$ is not spherically complete.  Following \cite{Wang}, fix a
maximal immediate extension $\widetilde K/K$, set
$\widetilde V=\mathcal O_{\widetilde K}$, and put
\[
  B=\widetilde V[[X]].
\]
The field $\widetilde K$ is spherically complete and
$v(\widetilde K^\times)=\mathbb R$. Thus $B$ is a GCD domain by
\cite[Theorem~1.1(i)]{Wang}.

We record the precise consequences of the construction in \cite{Wang} that are needed below, replacing the variable $T$ by $X$.

\begin{proposition}\cite[Section 4 and 5]{Wang}\label{prop:Wang-data}
There are nonzero elements $a_i\in\mathfrak m_V$, real numbers
$\alpha_i=v(a_i)>0$, and a nonzero $c\in\mathfrak m_V$ such that
\[
  \alpha_1>\alpha_2>\cdots,
  \qquad
  \alpha_i\longrightarrow0,
  \qquad
  v(c)=\sum_{i\geq1}\alpha_i.
\]
Set
\[
  \xi_i=a_i-X,
  \qquad
  P_n=\prod_{i=1}^n\xi_i.
\]
There are $F\in B$ and nonzero $f,g\in A$ with the following
properties:
\begin{enumerate}[label=\textup{(\roman*)}]
  \item $F(0)=c$ and $P_n\mid F$ in $B$ for every $n$;
  \item $P_n\mid f,g$ in $A$ for every $n$, and $F$ is a greatest
  common divisor of $f$ and $g$ in $B$;
  \item the ideal
  \[
    E=(f)A\cap(g)A
  \]
  is not finitely generated over $A$;
  \item there is a strictly increasing sequence of positive integers
  $(N_i)$ and an element $c_1\in K^\times$ such that every $N_i$ is
  prime to the residue characteristic when it is positive and
  \begin{equation}\label{eq:Wang-product}
    f(X)=c_1\prod_{i\geq1}
    \left(1-\left(\frac{X}{a_i}\right)^{N_i}\right)
  \end{equation}
  in $K[[X]]$.
\end{enumerate}
\end{proposition}

\begin{proof}
The valuation data are given by
\cite[Construction Lemma~4.2]{Wang}.  The remaining assertions follow
from \cite[Propositions~5.1--5.3 and the proof of
Theorem~4.1(i)]{Wang}; the normalized product
\eqref{eq:Wang-product} appears in the proof of
\cite[Proposition~5.2]{Wang}.
\end{proof}

We shall use evaluation at elements of positive value.

\begin{lemma}\label{lem:evaluation}
Let $C$ be the valuation ring of a complete nonarchimedean valued
field, and let $0\neq a\in\mathfrak m_C$.  Evaluation at $a$ defines a
surjective homomorphism
\[
  \operatorname{ev}_a:C[[X]]\longrightarrow C
\]
whose kernel is $(X-a)$.  Moreover, if $H\in X^NC[[X]]$, then
\begin{equation}\label{eq:evaluation-continuity}
  v\bigl(H(a)\bigr)\geq Nv(a).
\end{equation}
In particular, evaluation is continuous from the $X$-adic topology to
the valuation topology, and $X-a$ is prime in $C[[X]]$.
\end{lemma}

\begin{proof}
For $H(X)=\sum_{n\geq0}h_nX^n$, the series
$H(a)=\sum_{n\geq0}h_na^n$ converges in $C$.  Passage to the limit from
polynomial truncations shows that evaluation is a homomorphism.  The estimate
\eqref{eq:evaluation-continuity} is immediate.  The coefficient of
$X^r$ in
\[
  Q(X)=\sum_{n\geq1}h_n\sum_{j=0}^{n-1}X^{n-1-j}a^j
\]
is a convergent series in $C$, so $Q\in C[[X]]$.  The identity
\[
  H(X)-H(a)=(X-a)Q(X)
\]
shows that the kernel of evaluation is $(X-a)$.  Surjectivity follows
from the constant series.  The quotient is isomorphic to the domain
$C$, and hence $X-a$ is prime.
\end{proof}

For $n\geq1$, define
\begin{equation}\label{eq:Yn}
  Y_n=\frac{F}{P_n}\in B.
\end{equation}
Then
\begin{equation}\label{eq:Y-recursion}
  Y_n=\xi_{n+1}Y_{n+1}
\end{equation}
and, since $P_n(0)=\prod_{i=1}^na_i$,
\begin{equation}\label{eq:Y-constant-values}
  v\bigl(Y_n(0)\bigr)
  =v(c)-\sum_{i=1}^n\alpha_i
  =\sum_{i>n}\alpha_i
  \longrightarrow0.
\end{equation}

\begin{lemma}\label{lem:simple-linear-factors}
For every $n\geq1$,
\[
  \xi_n^2\nmid F
\]
in $B$.  Consequently,
\begin{equation}\label{eq:Ym-nonzero-evaluation}
  Y_m(a_n)\neq0
  \qquad(m\geq n).
\end{equation}
\end{lemma}

\begin{proof}
For $m\geq n$, let
\[
  f_m=c_1\prod_{i=1}^m
  \left(1-\left(\frac{X}{a_i}\right)^{N_i}\right)\in A.
\]
These are the partial products used in the construction
\cite[Proposition~5.2]{Wang}; they converge $X$-adically to $f$, and,
for each fixed $n$, one has $\xi_n\mid f_m$ for all sufficiently large $m$.
Since $\xi_n\bmod X=a_n\neq0$, division by $\xi_n$ on $(\xi_n)$ is
$X$-adically continuous by \cite[Lemma~2.3(i)]{Wang}.  Hence
$f_m/\xi_n$ converges $X$-adically to $f/\xi_n$.  Applying the
continuous evaluation map of Lemma~\ref{lem:evaluation} and using
\[
  \frac{1-(X/a_n)^{N_n}}{a_n-X}
  =\frac1{a_n^{N_n}}
    \sum_{r=0}^{N_n-1}a_n^{N_n-1-r}X^r,
\]
we obtain
\begin{equation}\label{eq:simple-factor-evaluation}
  \left(\frac{f}{\xi_n}\right)(a_n)
  =\frac{c_1N_n}{a_n}
   \prod_{\substack{j\geq1\\ j\neq n}}
   \left(1-\left(\frac{a_n}{a_j}\right)^{N_j}\right).
\end{equation}
This element is nonzero.  The finitely many factors with $j<n$ are
nonzero because
$v(a_n/a_j)=\alpha_n-\alpha_j<0$.  For $j>n$,
\[
  v\!\left(\left(\frac{a_n}{a_j}\right)^{N_j}\right)
  =N_j(\alpha_n-\alpha_j)\longrightarrow+\infty,
\]
so the tail product and its inverse converge.  By construction,
$N_n$ is a unit of $V$; this is automatic in residue characteristic
zero and follows from the prime-to-residue-characteristic condition in
positive residue characteristic.  Thus the right-hand side of
\eqref{eq:simple-factor-evaluation} is nonzero.  Hence
$\xi_n^2\nmid f$, and therefore $\xi_n^2\nmid F$ because $F\mid f$ in
$B$.

Let $m\geq n$.  If $Y_m(a_n)=0$, Lemma~\ref{lem:evaluation} gives
$\xi_n\mid Y_m$.  Since $P_m$ contains the factor $\xi_n$, the equality
$F=P_mY_m$ would imply $\xi_n^2\mid F$, a contradiction.
\end{proof}

The next construction produces the coprime quotients required by
Lemma~\ref{lem:flat-descent}.

\begin{proposition}\label{prop:coprime-tail}
There is an element $Z\in B$ such that
\begin{equation}\label{eq:coprime-tail}
  \gcd_B(Y_1,Z)=1
\end{equation}
and, with
\[
  s=\frac{fg}{F}\in B,
\]
one has
\begin{equation}\label{eq:tail-elements-in-E}
  sY_1,sZ\in E.
\end{equation}
\end{proposition}

\begin{proof}
Choose integers
\[
  0=k_1<k_2<k_3<\cdots
\]
recursively so that, whenever $1\leq n<m$,
\begin{equation}\label{eq:k-choice}
  (k_m-k_n)\alpha_n+v\bigl(Y_m(a_n)\bigr)
  >v\bigl(Y_n(a_n)\bigr)+m.
\end{equation}
At the $m$-th stage there are only finitely many inequalities, and
$\alpha_n>0$, so such an integer $k_m$ exists.

For $n\geq1$, put
\begin{equation}\label{eq:Zn}
  Z_n=\sum_{m\geq n}X^{k_m-k_n}Y_m\in B.
\end{equation}
The series converges $X$-adically and satisfies
\begin{equation}\label{eq:Z-recursion}
  Z_n=Y_n+X^{k_{n+1}-k_n}Z_{n+1}.
\end{equation}
By Lemma~\ref{lem:evaluation} evaluation at $a_n$ is continuous for
the $X$-adic topology.  Hence
\[
  Z_n(a_n)=\sum_{m\geq n}a_n^{k_m-k_n}Y_m(a_n).
\]
The first term is nonzero by
\eqref{eq:Ym-nonzero-evaluation}, while
\eqref{eq:k-choice} gives, for $m>n$,
\[
  v\bigl(a_n^{k_m-k_n}Y_m(a_n)\bigr)
  >v\bigl(Y_n(a_n)\bigr)+m.
\]
Thus the sum of the later terms has valuation strictly larger than
$v(Y_n(a_n))$, and the first term cannot be cancelled.  Therefore
\begin{equation}\label{eq:Zn-not-divisible}
  Z_n(a_n)\neq0,
  \qquad
  \xi_n\nmid Z_n.
\end{equation}

We prove that $Y_1$ and $Z_1$ are relatively prime.  Let $d\in B$ be a
common divisor.  Since $Y_1(0)\neq0$, one has $X\nmid d$.  The element
$X$ is prime because $B/XB\cong\widetilde V$ is a domain, so
$\gcd_B(d,X)=1$.  We claim inductively that
\begin{equation}\label{eq:d-divides-all}
  d\mid Y_n
  \quad\text{and}\quad
  d\mid Z_n
  \qquad(n\geq1).
\end{equation}
Suppose this holds for $n$.  Equation \eqref{eq:Z-recursion} gives
\[
  X^{k_{n+1}-k_n}Z_{n+1}\in dB\cap X^{k_{n+1}-k_n}B.
\]
Since $d$ is relatively prime to $X^{k_{n+1}-k_n}$ and relatively
prime elements $a,b$ of a GCD domain satisfy $aB\cap bB=abB$,
cancellation gives $d\mid Z_{n+1}$.  By
\eqref{eq:Zn-not-divisible}, $\xi_{n+1}\nmid Z_{n+1}$, and hence
$\xi_{n+1}\nmid d$.  Since $\xi_{n+1}$ is prime by
Lemma~\ref{lem:evaluation}, $d$ and $\xi_{n+1}$ are relatively prime.
Now
\[
  Y_n=\xi_{n+1}Y_{n+1}\in dB\cap\xi_{n+1}B
  =d\xi_{n+1}B,
\]
and cancellation gives $d\mid Y_{n+1}$.  This proves
\eqref{eq:d-divides-all}.

Write $Y_n=de_n$ with $e_n\in B$.  Since $Y_n(0)\neq0$, both
$d(0)$ and $e_n(0)$ are nonzero.  Taking constant terms gives
\[
  0\leq v\bigl(d(0)\bigr)\leq v\bigl(Y_n(0)\bigr)
  \qquad(n\geq1).
\]
By \eqref{eq:Y-constant-values}, $v(d(0))=0$.  Hence $d(0)$ is a unit
of $\widetilde V$, so $d$ is a unit of $B$.  This proves
\eqref{eq:coprime-tail} with $Z=Z_1$.

It remains to prove \eqref{eq:tail-elements-in-E}.  Since $P_m$
divides both $f$ and $g$ in $A$,
\[
  sY_1=\frac{fg}{P_1}\in E.
\]
Moreover,
\begin{align*}
  sZ
  &=\sum_{m\geq1}X^{k_m}\frac{fg}{P_m}\\
  &=f\left(\sum_{m\geq1}X^{k_m}\frac{g}{P_m}\right)
   =g\left(\sum_{m\geq1}X^{k_m}\frac{f}{P_m}\right).
\end{align*}
The two series in parentheses converge $X$-adically in $A$.
Consequently, $sZ\in(f)A\cap(g)A=E$.
\end{proof}

\begin{theorem}\label{thm:full-real-nonspherical}
Let $K$ be a complete nonarchimedean valued field with
$v(K^\times)=\mathbb R$, let $V=\mathcal O_K$, and assume that $K$ is
not spherically complete.  Set $A=V[[X]]$.  Then there are nonzero
$f,g\in A$ such that $E:=(f)\cap(g)$ is not flat over $A$.  Moreover,
\[
  \fd_A\bigl(A/(f,g)\bigr)\geq3,
  \qquad
  \wdim A\geq3.
\]
\end{theorem}

\begin{proof}
We take $f,g\in A$ to be the elements in Proposition~\ref{prop:Wang-data}, so that $E=(f)\cap (g)$ is not finitely generated over $A$. We will construct elements $y,z,s\in B$ such that the data $A,B,E,s,y,z$ satisfies the conditions of the Flat Descent criterion (Lemma~\ref{lem:flat-descent}), and it will follow that $E$ is not flat.

Write $f=Ff_0$ and $g=Fg_0$ in $B$. Since $F$ is a greatest common divisor,  $f_0$ and $g_0$ are relatively prime in $B$, so $f_0B\cap g_0B=f_0g_0B$. So writing $s:=Ff_0g_0=\frac{fg}{F}$, we see that
\[
  (f)B\cap(g)B
  =F(f_0B\cap g_0B)
  =Ff_0g_0B
  sB.
\]
It follows that $E\subseteq sB$.

The natural embedding
$\Frac(A)\subseteq K((X))$ gives
\[
  \Frac(A)\cap B
  \subseteq K((X))\cap\widetilde V[[X]].
\]
A Laurent series in the last intersection has no negative powers and
all its coefficients lie in
$K\cap\widetilde V=V$.  Thus $\Frac(A)\cap B\subseteq A$, and we see immediately that
\begin{equation}\label{eq:full-real-contraction}
  \Frac(A)\cap B=A.
\end{equation}
Let $y:=Y_1=\frac{F}{P_1}$, and take $z\in B$ to be the element $Z$ given by Proposition~\ref{prop:coprime-tail}, and it follows that this data satisfies the conditions of Lemma~\ref{lem:flat-descent} as required.\end{proof}

\section{Completion invariance}\label{sec:completion}

\noindent In this section, we will complete the proof of our main theorem. In light of Theorem \ref{thm:full-real-nonspherical}, the only remaining case to consider is when $V$ is non-complete, with value group $\mathbb{R}$. The valuation on $\widehat K$ is again denoted by $v$.

Again, let $\widehat K$ be the metric completion of $K=Frac(V)$ and
$\widehat V=\mathcal O_{\widehat K}$ for its valuation ring. It follows from Lemma \ref{lem: maximal completion} that $\widehat{K}$ is spherically complete if and only if $V$ is almost maximal.

\begin{lemma}\label{lem:completion-congruences}
The extension $\widehat K/K$ is immediate.  Moreover, for every
$0\neq a\in V$, the natural map
\begin{equation}\label{eq:completion-congruences}
  V/aV\longrightarrow\widehat V/a\widehat V
\end{equation}
is an isomorphism.
\end{lemma}

\begin{proof}
Let $0\neq x\in\widehat K$.  By density, there is $y\in K$ such that
$v(x-y)>v(x)$, and then $v(y)=v(x)$.  Thus completion does not enlarge
the value group.  If $x\in\widehat V$, choose $y\in K$ with
$v(x-y)>0$.  Then $y\in V$ and $x,y$ have the same residue class.
Hence the residue fields also agree, so $\widehat K/K$ is immediate.

For surjectivity in \eqref{eq:completion-congruences}, let
$x\in\widehat V$.  Choose $y\in K$ with $v(x-y)>v(a)$.  Then $y\in V$
and $x-y\in a\widehat V$.  For injectivity, if
$y\in V\cap a\widehat V$, then $y/a\in K\cap\widehat V=V$, and hence
$y\in aV$.
\end{proof}

We will now consider Laurent series rings over $V$ and $\widehat{V}$. Put
\[
  R=V((X)),\qquad S=\widehat V((X)).
\]

\begin{proposition}\label{prop:Laurent-completion-flat}
For every $0\neq t\in R$, the natural map
\begin{equation}\label{eq:principal-quotient-completion}
  R/tR\longrightarrow S/tS
\end{equation}
is an isomorphism.  The inclusion $R\to S$ is faithfully flat.
\end{proposition}

\begin{proof}
Since $X$ is a unit in both rings, after multiplying $t$ by a power of
$X$ we may write
\[
  t=\sum_{i\geq0}t_iX^i\in V[[X]],\qquad t_0\neq0.
\]
Let $s=\sum_{n\geq N}s_nX^n\in\widehat R$.  We construct
$b_n\in V$ and $y_n\in\widehat V$ recursively so that, with
$b=\sum_{n\geq N}b_nX^n$ and $y=\sum_{n\geq N}y_nX^n$, one has
$s=b+ty$.  Set $y_j=0$ for $j<N$.  Once $y_N,\ldots,y_{n-1}$ have
been chosen, put
\[
  r_n=s_n-\sum_{i=1}^{n-N}t_i y_{n-i}\in\widehat V.
\]
Lemma~\ref{lem:completion-congruences}, applied to $a=t_0$, provides
$b_n\in V$ such that $r_n-b_n\in t_0\widehat V$.  Define
\[
  y_n=\frac{r_n-b_n}{t_0}\in\widehat V.
\]
The coefficient recursion gives $s=b+ty$.  Hence
$S=R+tS$, proving surjectivity of
\eqref{eq:principal-quotient-completion}.

For injectivity, suppose that $b=ty$ with $b\in R$ and
$y\in S$.  If $y\neq0$, write
$y=\sum_{n\geq N}y_nX^n$ with $y_N\neq0$.  Comparing coefficients and
arguing inductively gives
\[
  t_0y_n=b_n-\sum_{i=1}^{n-N}t_i y_{n-i}\in K
  \qquad(n\geq N).
\]
Thus every $y_n$ lies in $K\cap\widehat V=V$, so $y\in R$.  Therefore
\begin{equation}\label{eq:principal-contraction-completion}
  R\cap tS=tR,
\end{equation}
and \eqref{eq:principal-quotient-completion} is an isomorphism.

Let $Q=S/R$.  The identities
$S=R+tS$ and
\eqref{eq:principal-contraction-completion} show that multiplication by
every $0\neq t\in R$ is bijective on $Q$.  Hence $Q$ is naturally a
vector space over $\Frac(R)$, and is therefore flat over $R$.  The
exact sequence
\[
  0\longrightarrow R\longrightarrow S
  \longrightarrow Q\longrightarrow0
\]
then shows that $S$ is flat over $R$.  Since $Q$ is flat,
tensoring this sequence with any $R$-module $M$ yields an injection
\[
  M\longrightarrow M\otimes_R S.
\]
Thus the extension is faithfully flat.
\end{proof}

\begin{lemma}\label{lem:finitely-generated-ideal-descent}
Every nonzero element of $S$ is associated in $S$ to
an element of $R$.  Consequently, every finitely generated ideal of
$S$ is extended from a finitely generated ideal of $R$.
\end{lemma}

\begin{proof}
Let $0\neq c\in S$.  Multiplication by a power of $X$, which is
a unit in both rings, reduces us to
\[
  c=\sum_{n\geq0}c_nX^n\in\widehat V[[X]],\qquad c_0\neq0.
\]
Since $\widehat K/K$ is immediate, choose $0\neq b_0\in V$ with
$v(b_0)=v(c_0)$.  Then $b_0\widehat V=c_0\widehat V$ and
$u_0=b_0/c_0$ is a unit of $\widehat V$.

We construct $b_n\in V$ and $u_n\in\widehat V$ recursively so that
\[
  \left(\sum_{i\geq0}c_iX^i\right)
  \left(\sum_{j\geq0}u_jX^j\right)
  =\sum_{n\geq0}b_nX^n.
\]
Suppose $u_0,\ldots,u_{n-1}$ have been chosen and set
\[
  r_n=\sum_{i=1}^n c_i u_{n-i}\in\widehat V.
\]
By Lemma~\ref{lem:completion-congruences}, choose $b_n\in V$ such that
$b_n-r_n\in b_0\widehat V=c_0\widehat V$, and put
\[
  u_n=\frac{b_n-r_n}{c_0}\in\widehat V.
\]
Then, for
$u=\sum_{n\geq0}u_nX^n$ and $b=\sum_{n\geq0}b_nX^n$, one has
$cu=b$.  Since $u_0$ is a unit, $u\in S^\times$.  Hence $c$
and $b\in R$ are associated in $S$.

If $J=(c_1,\ldots,c_m)S$, choose $b_i\in R$ associated to
$c_i$.  Then
\[
  J=(b_1,\ldots,b_m)S,
\]
which proves the final assertion.
\end{proof}

\begin{theorem}\label{thm:completion-invariance}
$\wdim V((X))=\wdim\widehat V((X))$.
\end{theorem}

\begin{proof}
Again, write $R=V((X))$ and $S=\widehat V((X))$.  Since $S$ is faithfully
flat over $R$ by Proposition \ref{prop:Laurent-completion-flat}, flat base change gives, for $R$-modules $M,N$,
\[
  \Tor_n^R(M,N)\otimes_RS
  \cong
  \Tor_n^S(M\otimes_RS,N\otimes_RS).
\]
Faithfulness therefore implies
\begin{equation}\label{eq:first-wdim-inequality}
  \wdim R\leq\wdim S.
\end{equation}

Suppose first that $d=\wdim R<\infty$.  Let $M$ be an $S$-module and
let $J\subseteq S$ be a finitely generated ideal.  By
Lemma~\ref{lem:finitely-generated-ideal-descent}, there is a finitely
generated ideal $I\subseteq R$ such that $J=IS$.  A free resolution of
$R/I$, after tensoring with the flat $R$-algebra $S$, gives
\begin{equation}\label{eq:Tor-change-completion}
  \Tor_n^S(S/J,M)\cong\Tor_n^R(R/I,M)
  \qquad(n\geq0),
\end{equation}
where $M$ on the right is viewed as an $R$-module.  Thus
$\Tor_{d+1}^S(S/J,M)=0$ for every finitely generated ideal $J$ of $S$.
Let $L$ be the $d$-th syzygy of $M$ in a free resolution, with the
convention $L=M$ when $d=0$.  Dimension shifting gives
\[
  \Tor_1^S(S/J,L)=0
\]
for every finitely generated $J$.  The equational criterion for
flatness now shows that $L$ is flat.  Hence $\fd_SM\leq d$, and
$\wdim S\leq\wdim R$.  Together with
\eqref{eq:first-wdim-inequality}, this proves
that $\wdim R=\wdim S$ when $\wdim R$ is finite.  If
$\wdim R=\infty$, \eqref{eq:first-wdim-inequality} forces
$\wdim S=\infty$, so the same equality holds.\qedhere

\end{proof}

\begin{corollary}\label{cor: value group R}
For any valuation ring $V$ of rank 1, with value group $\mathbb{R}$, $\wdim V[[X]]=2$ if and only if $V$ is almost maximal, otherwise $\wdim V[[X]] = \wdim V((X))>2$.
\end{corollary}

\begin{proof}

If $V$ is almost maximal, then $V[[X]]$ is coherent by Theorem \ref{thm: coherence}, so $\wdim V[[X]]=2$ by \cite[Lemma 1]{JondrupSmall}.

Conversely, if $V$ is not almost maximal, then $\widehat{K}$ is spherically complete by Lemma \ref{lem: maximal completion}, so it follows from Theorem \ref{thm:full-real-nonspherical} that $\wdim\widehat{V}[[X]]\geq 3$. But since $\wdim\widehat{V}=1$, this means that $\wdim \widehat{V}[[X]]>\wdim\widehat{V}+1$, and hence $\wdim\widehat{V}[[X]]=\wdim\widehat{V}((X))>2$ by Proposition \ref{prop:known-localization}.\\

\noindent But $\wdim\widehat{V}((X))=\wdim V((X))>2$ by Theorem \ref{thm:completion-invariance}, and applying Proposition \ref{prop:known-localization} again gives $$\wdim V[[X]]=\max\{\wdim V+1,\wdim V((X))\}=\max\{2,\wdim V((X))\}=\wdim V((X))\qedhere$$
\end{proof}

\noindent Combining all our results so far allow us to finally complete the proof of our main result.\\

\noindent\emph{Proof of Theorem \ref{main_theorem}.} Using Corollary \ref{cor: conditions for coherence}, we know $i\iff iii$, and $i\implies ii$ is given by \cite[Lemma 1.1]{JondrupSmall}, so it remains to prove $ii\implies iii$. We will assume $iii$ is false, i.e. either $V$ is not almost maximal or $\Gamma\neq 0,\mathbb{R},\mathbb{Z}$, and show that $\wdim V[[X]]=\wdim V((X))>2$, which will complete the proof.\\

\noindent If $\Gamma\neq 0,\mathbb{R}$ and $\Gamma\not\cong\mathbb{Z}$, then $\Gamma$ is a proper dense subgroup of $\mathbb{R}$, so it follows from Theorem \ref{thm:main'} that $\wdim V[[X]]>2$. Since $\wdim V=1$, this means that $\wdim V[[X]]>\wdim V+1$, and hence $\wdim V[[X]]=\wdim V((X))>2$ by Proposition \ref{prop:known-localization}.\\

\noindent On the other hand, if $V$ is not almost maximal, then $V$ is not discrete, so we may assume that $\Gamma=\mathbb{R}$, and hence $\wdim V[[X]]=\wdim V((X))>2$ by Theorem \ref{cor: value group R}.\qed

\begingroup
\sloppy
\hbadness=10000

\endgroup

\end{document}